\documentclass[12pt,oneside,a4paper,abstracton]{scrartcl}
\usepackage{graphicx}
\usepackage{amsmath}
\usepackage{amssymb}
\usepackage[table]{xcolor}
\usepackage[colorlinks,linkcolor=blue,urlcolor=blue,citecolor=blue]{hyperref}

\usepackage{caption,subcaption} 

\newcommand{\cR}{\ensuremath{\mathcal{R}}}
\newcommand\R{\ensuremath{\mathbb{R}}}
\newcommand{\Rnn}{\ensuremath{\R^{n\times n}}}
\newcommand{\Rnk}{\ensuremath{\R^{n\times k}}}
\newcommand{\Rkk}{\ensuremath{\R^{k\times k}}}
\newcommand{\Rnm}{\ensuremath{\R^{n\times m}}}
\newcommand{\Rqn}{\ensuremath{\R^{q\times n}}}

\newcommand{\diag}{\ensuremath{\operatorname{diag}}}

\newcommand{\cM}{\ensuremath{\mathcal{M}}}
\newcommand{\cT}{\ensuremath{\mathcal{T}}}

\newcommand{\bm}[1]{\mathbf{#1}}
\newcommand{\bX}{\ensuremath{\bm{X}}}
\newcommand{\bY}{\ensuremath{\bm{Y}}}
\newcommand{\ba}{\ensuremath{\bm{a}}}
\newcommand{\bb}{\ensuremath{\bm{b}}}
\newcommand{\bg}{\ensuremath{\bm{g}}}
\newcommand{\hA}{\ensuremath{\hat{A}}}
\newcommand{\hD}{\ensuremath{\hat{D}}}
\newcommand{\hL}{\ensuremath{\hat{L}}}

\newcommand{\ttt}{\ensuremath{\tilde{t}}}
\newcommand{\tA}{\ensuremath{\tilde{A}}}
\newcommand{\tD}{\ensuremath{\tilde{D}}}
\newcommand{\tS}{\ensuremath{\tilde{S}}}
\newcommand{\tW}{\ensuremath{\tilde{W}}}
\newcommand{\tX}{\ensuremath{\tilde{X}}}

\newcommand{\lyapack}{\textsf{LyaPack}}
\newcounter{mymac@matlab}
\newcommand{\matlab}{MATLAB%
  \ifnum\value{mymac@matlab}<1%
  \textsuperscript{\textregistered}%
  \setcounter{mymac@matlab}{1}%
  \fi
}

\title{Peer Methods for the Solution of Large-Scale Differential Matrix
  Equations}%
\author{Peter Benner \and Norman Lang}
\date{\today}%
\begin{document}

\maketitle

\noindent
{\bfseries Author's addresses:}\\[2ex]
Peter Benner \\
Max Planck Institute for Dynamics of Complex Technical Systems,\\ 
Computational Methods in Systems and Control Theory,\\
D-39106 Magdeburg, Germany,\\
and\\
Technische Universit\"at Chemnitz,\\ 
Faculty of Mathematics, Mathematics in Industry and Technology\\
D-09126 Chemnitz, Germany\\
\url{benner@mpi-magdeburg.mpg.de}\\[2ex]
Norman Lang\\
Technische Universit\"at Chemnitz,\\
Faculty of Mathematics, Mathematics in Industry and Technology\\
D-09126 Chemnitz, Germany\\
\url{norman.lang@mathematik.tu-chemnitz.de}

\begin{abstract}
  We consider the application of implicit and linearly implicit
  (Rosenbrock-type) peer methods to matrix-valued ordinary differential
  equations. In particular the differential Riccati equation
  (DRE) is investigated. For the Rosenbrock-type schemes, a reformulation
  capable of avoiding a number of Jacobian applications is developed that, in
  the autonomous case, reduces the computational complexity of the
  algorithms. Dealing with large-scale problems, an efficient implementation
  based on low-rank symmetric indefinite factorizations is presented. The
  performance of both peer approaches up to order 4 is compared to existing
  implicit time integration schemes for matrix-valued differential equations.
\end{abstract}
\section{Introduction}\label{intro}
Differential matrix equations are of major importance in many fields like
optimal control and model reduction of linear dynamical systems, see,
e.g.,~\cite{Loc01,AboFIetal03} and \cite{ShoSV83,VerK83}, respectively. In that
context, the most common differential matrix equations are the differential
Riccati and Lyapunov equations (DREs/DLEs), where the latter can be considered as a
special case of the Riccati equation. Therefore, as an illustrating example, in
this article, we consider the numerical solution of the time-varying DRE
\begin{align}
  \begin{aligned}
    \dot{X}(t)&=A(t)^{T}X(t)+X(t)A(t)-X(t)S(t)X(t)+W(t)=:\cR(t,X),\\
    X(t_0)&=X_{0},
  \end{aligned}\label{eq:DRE}
\end{align}
where $t\in[t_{0},t_{f}]$, $X(t)\in\Rnn$ is the sought for solution to
Equation~(\ref{eq:DRE}) and $A(t)$, $W(t)$, $S(t)\in\Rnn$ are given
matrix-valued functions and the matrix $X_{0}\in\Rnn$ denotes the initial value
with $n$ being the problem dimension. The differential Lyapunov equation results
if we set \(S(t)\equiv 0\). Provided that the matrices $A,W,S$ are
piecewise continuous and locally bounded, from~\cite[Theorem 4.1.6]{AboFIetal03}
we have that a solution to Equation~(\ref{eq:DRE}) exists and further is
unique.

The DRE is one of the most deeply studied nonlinear matrix differential
equations due to its importance in optimal control, optimal filtering,
$\mathbf{H}_\infty$-control of linear-time varying systems, differential games
and many more (see, e.g.,~\cite{AboFIetal03,IchK99,Jac93,PetUS00}). Over the
last four decades many solution strategies have been presented, see, e.g.,
\cite{DavM73,Lai76,Lau82,HarWA83,KenL85,Meh91}. Most of these methods are only
applicable to small-scale systems, i.e., systems with a rather small number $n$
of unknowns. Others, suitable for the application to large-scale problems have
to deal with numerical difficulties like, e.g., instability, see~\cite[Section
4.1]{Men07} for a detailed overview. Due to the fact that in many control
applications fast and slow modes are present, the DRE~(\ref{eq:DRE}) is usually
fairly stiff. For that reason, the numerical solution based on matrix versions
of classical implicit time integration schemes, such as the BDF, Rosenbrock
methods, and the Midpoint and Trapezoidal
rules~\cite{Die92,BenM04,Men07,BenM13,BenM18} has become a popular tool for the
solution of~(\ref{eq:DRE}). Recently, also splitting
methods~\cite{HanS14,Sti15,Sti17} and a structure preserving solution method for
large-scale DREs~\cite{KosM17}, using Krylov subspace methods, were
proposed. Note that ``unrolling'' the matrix differential equation into a
standard (vector-valued) ordinary differential equation (ODE) is usually
infeasible due to the resulting complexity - the ODE would then be posed in
\(\R^{n^{2}\times n^{2}}\).

In the field of implicit time integration methods, linear multistep and one-step
methods, such as e.g., the BDF and the Rosenbrock methods, respectively, have
been known for many decades and in addition have proven their effectiveness over
the years for a wide range of problems. These two traditional classes of time
integration methods have been studied separately until recently. As a unifying
framework for stability, consistency and convergence analysis for a wide variety
of methods, also containing the aforementioned classes, in~\cite{But66} the
general linear methods (GLMs) were introduced. Detailed explanations on GLMs are
given in the surveys~\cite{But85,But96}. Most of the classical methods contain a
number of solution variables $X_{k+1,j}$ and in addition compute a separate
number of auxiliary variables related to function evaluations
$F(\ttt_{k},\tX_{k})$ at $t_{k}\leq \ttt_{k}\leq t_{k+1}$,
$\tX_{k}\approx X(\ttt_{k})$ that are designated to improve the accuracy and
stability properties of the approximate solutions. In particular, usually only
one solution variable for the approximation of the solution in each time
interval is employed. Moreover, for different time intervals, as, e.g., for
variable time step sizes, these solution variables may have distinguished
accuracy and stability properties. Due to that e.g., the Rosenbrock methods
often suffer from order reduction. Now, the idea of the so-called \emph{peer
  methods} is to define an integration scheme that only contains peer variables,
each representing an approximation to the exact solution of~\eqref{eq:DRE} at
different time locations that share the same accuracy and stability properties.

For all the above mentioned implicit solution methods, including the peer
methods to be presented, it turns out that the main ingredient for the solution
of the DRE~(\ref{eq:DRE}) is to solve a number of either algebraic Riccati or
Lyapunov equations (AREs/ALEs) in each time step. Dealing with large-scale
systems, the simple application of the implicit integration methods leads to an
enormous computational effort and storage amount in the sense that the solution
to the DRE~(\ref{eq:DRE}) is a dense square matrix of dimension $n$ being
computed at each point of the discrete time set. However, in practice it is
often observed that the singular values of the solution of the ALEs occurring in
the innermost iteration of the solution methods decay rapidly to zero, see,
e.g.,~\cite{AntSZ02,Gra04,Pen00a,TruV07}. 
Thus, the solution is of low numerical rank, meaning it can be well approximated
by products of low-rank matrices. Based on this observation, modern and
efficient algorithms rely on low-rank based solution
algorithms. In~\cite{Men07,BenM13,BenM18} classical implicit time integration
methods, originally developed for standard scalar and vector-valued ODEs,
exploiting the low-rank phenomena were presented. Therein a factorization
$X=ZZ^{T}$ with $Z\in\Rnk$, $k\ll n$, is employed in order to efficiently solve
large-scale DREs. This decomposition will be referred to as a low-rank
Cholesky-type factorization (LRCF) in the remainder. In~\cite{LanMS14,LanMS15},
it has be shown that for integration methods of order $\geq 2$, the right hand
sides of the ALEs to be solved become indefinite and thus the LRCF involves
complex data and therefore requires complex arithmetic and storage.  Moreover,
therein a low-rank symmetric indefinite factorization (LRSIF) of the form
$X=LDL^{T}$, $L\in\Rnk,~D\in\Rkk$, $k\ll n$, of the solution, was introduced in
order to avoid complex data.

The paper is organized as follows. In Section~\ref{sec:peer}, the implicit and
linearly implicit Ro\-sen\-brock-type peer methods are introduced for the
application to the matrix-valued differential Riccati equation. Efficient
numerical algorithms based on the low-rank symmetric indefinite factorization
for both peer approaches are presented in Section~\ref{sec:low-rank}. In
Section~\ref{sec:numerics}, the performance of the new peer solution methods up
to order 4 is compared to the existing implicit integration schemes of similar
orders. A conclusion is given in Section~\ref{sec:conc}.
\section{Peer Methods}\label{sec:peer}
The class of peer methods first appeared in~\cite{SchW04} in terms of linearly
implicit integration schemes with peer variables, suitable for parallel
computations by only using information from the previous time interval. A number
of specific peer schemes and applications are presented in,
e.g.,~\cite{PodWS05,PodWS06,SchWE05,SchWP05}. Further, for a recent detailed
overview see~\cite[Chapters 5,10]{StrWP12}.
\subsection{General Implicit Peer Methods}\label{sec:implPeer}
A general (one-step) implicit peer method, applied to the matrix-valued initial
value problem~\eqref{eq:DRE} reads
\begin{align}
  X_{k,i}=&\sum_{j=1}^{s}b_{i,j}X_{k-1,j}
  +\tau_{k}\sum_{j=1}^{i}g_{i,j}\cR(t_{k,j},X_{k,j}).\label{eq:peer}
\end{align}
Here, $s$ is the number of stages and
\begin{align}
  t_{k,j}=t_{k}+c_{j}\tau_{k},\label{eq:peertime}
\end{align}
where the variables $c_{j},~j=1,\dots,s$, with $c_{s}=1,~t_{k,s}=t_{k+1}$, define
the locations of the peer variables \(X_{k,i},~i=1,\dots,s\), for the time step
$t_{k}\rightarrow t_{k+1}$. In general, $c_{j}<0$ for some $j$ will also be
allowed. Furthermore, the peer variables $X_{k,i}$ represent the solution
approximations of~(\ref{eq:DRE}) at the time locations $t_{k,i}$, i.e,
$X_{k,i}\approx X(t_{k,i})=X(t_{k}+c_{i}\tau_{k})$. From $c_{s}=1$, the solution
$X_{k}$ at time $t_{k}$ is given by $X_{k-1,s}$. The variables \(b_{i,j}\) and
\(g_{i,j}\) are the determining coefficients of the method.

The convergence order of these methods is restricted to $s-1$. Thus,
additionally using function values from the previous time interval, two-step
peer methods of order $s$ can be constructed. Under special conditions even a
superconvergent subclass of the implicit peer methods with convergence order
$s+1$ can be found. Details on the convergence analysis are given in,
e.g.,~\cite{SolW17} and the references therein. The corresponding two-step scheme
becomes
\begin{align}
  X_{k,i}=&\sum_{j=1}^{s}b_{i,j}X_{k-1,j}+
            \tau_{k}\sum_{j=1}^{s}a_{i,j}\cR(t_{k-1,j},X_{k-1,j})
  +\tau_{k}\sum_{j=1}^{i}g_{i,j}\cR(t_{k,j},X_{k,j})\label{eq:peer2}
\end{align}
with additional coefficients \(a_{i,j}\).

Note that, given from the order conditions, the coefficients will in general
depend on the step size ratio $\tau_{k}/\tau_{k-1}$ of two consecutive time
steps. Moreover, the computation of the coefficients is based on highly
sophisticated optimization processes. Therefore, for details on the order
conditions and the computation of the associated coefficients, we refer
to~\cite{SolW17} and the references therein. Further, note that the
scheme~\eqref{eq:peer} can easily be recovered from~\eqref{eq:peer2} by setting
$a_{i,j}=0$ and therefore in the remainder the statements restrict to the more
general class~(\ref{eq:peer2}) of implicit (two-step) peer methods.

From the application of the peer scheme~\eqref{eq:peer2} to the DRE~\eqref{eq:DRE}
with $F(t_{k,i},X_{k,i})=\cR(t_{k,i},X_{k,i})$ one obtains
\begin{align}
  \tA_{k,i}^TX_{k,i}+X_{k,i}\tA_{k,i}
  -X_{k,i}\tS_{k,i}X_{k,i}+\tW_{k,i}=0,\qquad i=1,\dots,s,\label{eq:peerARE}
\end{align}
that in fact is an algebraic Riccati equation. Here, the coefficient matrices are
given by
\begin{align*}
  \tA_{k,i}&=\tau_{k}g_{i,i}A_{k,i}-\frac{1}{2}I,\quad
  \tS_{k,i}=\tau_{k}g_{i,i}S_{k,i},\\
  \tW_{k,i}&=\tau_{k}g_{i,i}W_{k,i}+\sum_{j=1}^{s}b_{i,j}X_{k-1,j}
  +\tau_{k}\sum_{j=1}^{s}a_{i,j}\cR(t_{k-1,j},X_{k-1,j})
  +\tau_{k}\sum_{j=1}^{i-1}g_{i,j}\cR(t_{k,j},X_{k,j}).
\end{align*}
Moreover, we have $A_{k,i}=A(t_{k,i})$, $W_{k,i}=W(t_{k,i})$ and
$S_{k,i}=S(t_{k,i})$ with $t_{k,i}$ from Equation~(\ref{eq:peertime}). Note
that, according to the number of peer variables to be computed, $s$ AREs have to
be solved at every time step $t_{k}\to t_{k+1}$ of the method.

In comparison, the BDF methods, as well as the Midpoint and Trapezoidal rules,
require the solution of only one ARE at every time step, see
e.g.,~\cite{Men07,LanMS15,Lan17}. That is, directly solving the occurring
algebraic Riccati equations, the expected computational effort of the peer
methods is in general $s$-times higher than that of the other implicit
methods. Still, from the fact that $s$ peer variables with the same accuracy and
stability properties are computed within every time interval, the peer methods
allow us to use larger step sizes in order to achieve a comparable accuracy. A
comparison and detailed investigation is given in
Section~\ref{sec:numerics}. Note that analogously to the DRE case, the peer
method can be applied to differential Lyapunov equations or any other
differential matrix equation. The application to differential Lyapunov equations
is presented in~\cite{Lan17}.

For solving the AREs, in general, any solution method suitable for sparse
large-scale problems can be applied. A detailed overview can, e.g., be found
in~\cite{BenS13,Sim16}. In this contribution, Newton's method is going to be
used in order to find a solution to the AREs~(\ref{eq:peerARE}) arising within
the peer scheme~(\ref{eq:peer2}). Following~\cite{Kle68,LanR95}, Newton's method
applied to the AREs~(\ref{eq:peerARE}) results in the solution of an algebraic
Lyapunov equation
\begin{align}
  {{}\hA_{k,i}^{(\ell)}}^TX_{k,i}^{(\ell)}+X_{k,i}^{(\ell)}\hA_{k,i}^{(\ell)}
  =-\tW_{k,i}-X_{k,i}^{(\ell-1)}\tS_{k,i}X_{k,i}^{(\ell-1)}
  \label{eq:ALE_implPeer}
\end{align}
with
\begin{align*}
  \hA_{k,i}^{(\ell)}&=\tA_{k,i}-\tS_{k,i}X_{k,i}^{(\ell-1)}.
\end{align*}
at each step $\ell$ of the Newton iteration and thus the solution of
\eqref{eq:DRE}, using the implicit peer scheme~\eqref{eq:peer2} boils down to
the solution of a sequence of ALEs at every time step of the integration scheme.
\subsection{Rosenbrock-Type Peer Methods}\label{sec:RosPeer}
\subsubsection{Standard Representation}\label{sec:stand-repr}
For the implicit peer methods applied to the DRE, a number of AREs has to be
solved.  In order to avoid the solution of these nonlinear matrix equations, we
also consider linearly implicit peer methods in terms of the two-step
Rosenbrock-type peer schemes
\begin{align}
  \begin{aligned}
    (I-\tau_{k}g_{i,i}J_{k})X_{k,i}=&\sum_{j=1}^{s}b_{i,j}X_{k-1,j}+
    \tau_{k}\sum_{j=1}^{s}a_{i,j}\left(F(t_{k-1,j},X_{k-1,j})-J_{k}X_{k-1,j}\right)\\
    &+\tau_{k}J_{k}\sum_{j=1}^{i-1}g_{i,j}X_{k,j},
  \end{aligned}\label{eq:peerRos}
\end{align}
introduced in~\cite{PodWS05}. As for the implicit schemes, here we consider
methods with\\ \mbox{$g_{1,1}=\dots=g_{s,s}=\gamma$}. For the comprehensive
derivation of coefficients $a_{i,j},~b_{i,j}$ and $g_{i,j}$ that result in
stable schemes~(\ref{eq:peerRos}), for arbitrary step size ratios, we refer
to~\cite[Section 3]{PodWS05}.  Expression $J_{k}$ denotes the Jacobian
represented by the Fr\'{e}chet derivative
\begin{align}
  J_{k}:=\frac{\partial \cR}{\partial X}(t_k,X_k): U \rightarrow
  (A_k-S_kX_k)^TU+U(A_k-S_kX_k).\label{eq:frechet}
\end{align}
of $F$ at $(t_{k},X_{k})$. Now, replacing the Jacobian $J_{k}$
in~(\ref{eq:peerRos}) by~(\ref{eq:frechet}), for the solution of the DRE, the
procedure reads
\begin{align}
  \begin{aligned}
    \tA_{k,i}^{T}X_{k,i}&+X_{k,i}\tA_{k,i}=-\tW_{k,i}\quad i=1,\dots,s,\\
    \tW_{k,i}=&\sum_{j=1}^{s}b_{i,j}X_{k-1,j}
    +\tau_{k}\sum_{j=1}^{s}a_{i,j}\left(\cR(t_{k-1,j},X_{k-1,j})
      -(\hA_{k}^{T}X_{k-1,j}+X_{k-1,j}\hA_{k})\right)\\
    &+\tau_{k}\sum_{j=1}^{i-1}g_{i,j}(\hA_{k}^{T}X_{k,j}+X_{k,j}\hA_{k})
  \end{aligned}\label{eq:ALE_RosPeer}
\end{align}
with the matrices $\hA_{k}=A_{k}-S_{k}X_{k}$ and
$\tA_{k,i}=\tau_{k}g_{i,i}\hA_{k}-\frac{1}{2}I$.

\subsubsection{Reformulation to avoid Jacobian applications}
\label{sec:jacob-avoid-reform}
The Rosenbrock-type peer scheme~(\ref{eq:peerRos}) involves the solution
of an ALE at each stage. The right hand sides of these ALEs
particularly require the application of the Jacobian
$J_{k}$ to the sums $\sum_{j=1}^{s}a_{i,j}X_{k-1,j}$ and
$\sum_{j=1}^{i}g_{i,j}X_{k,j}$ of the previous and current solution
approximations, respectively. In order to at least avoid the application of the
Jacobians to the sum of new variables $X_{k,j}$, a reformulation, similar to
what is standard for the classical Rosenbrock schemes, see, e.g.,~\cite[Chapter
IV.7]{HaiW02}, based on the variables
\begin{align}
  Y_{k,i}=\sum_{j=1}^{i}g_{i,j}X_{k,j},~i=1,\dots,s\Leftrightarrow
  \bY_{k}=(G\otimes I_{n})\bX_{k}\label{eq:auxvar_peer}
\end{align}
can be stated. Here,
$\bX_{k}=(X_{k,i})_{i=1}^{s},~\bY_{k}=(Y_{k,i})_{i=1}^{s}\in\R^{sn\times n}$ and
$\otimes$ denotes the Kronecker product. Provided that $g_{i,i}\neq 0,~\forall
i$, the lower triangular matrix $G=(g_{i,j})$ is non-singular and the original
variables $X_{k,i}$ can be recovered from the relation
\begin{align}
  \bX_{k}=(G^{-1}\otimes I_{n})\bY_{k}\Leftrightarrow
  X_{k,i}=\sum_{j=1}^{i}\bg_{i,j}Y_{k,j},~i=1,\dots,s\label{eq:origvar_peer}
\end{align}
where $G^{-1}=(\bg_{i,j})$ and $\bg_{i,i}=\frac{1}{g_{i,i}}$. Then,
from~(\ref{eq:origvar_peer}), we obtain
\begin{align}
  \begin{aligned}
    \sum_{j=1}^{s}\!a_{i,j}X_{k-1,j}, i=1,\dots,s
    \Leftrightarrow\quad
    \!&~((a_{i,j})\otimes I)\bX_{k}&\\
    \!=&~((a_{i,j})\otimes I)(G^{-1}\!\otimes I)\bY_{k}&\\
    \!=&~((a_{i,j})G^{-1}\otimes I)\bY_{k}
    &\!\!\!\!\!\Leftrightarrow\!\sum_{j=1}^{s}\!\ba_{i,j}Y_{k-1,j},i=1,\dots,s
  \end{aligned}\label{eq:modsumPeer}
\end{align}
with the coefficients
\begin{align}
  \begin{aligned}
    (\ba_{i,j})=(a_{i,j})G^{-1}.
  \end{aligned}\label{auxcoeff_peer}
\end{align}
and analogously, for the sum $\sum_{j=1}^{s}b_{i,j}X_{k-1,j}$, we have
\begin{align*}
  \sum_{j=1}^{s}b_{i,j}X_{k-1,j}=\sum_{j=1}^{s}\bb_{i,j}Y_{k-1,j},~i=1,\dots,s
\end{align*}
with $(\bb_{i,j})=(b_{i,j})G^{-1}$.

Now, inserting the auxiliary variables~(\ref{eq:auxvar_peer})
into~(\ref{eq:peerRos}) and dividing the result by $\tau_{k}$, the linearly
implicit scheme can be reformulated to
\begin{align}
  \begin{aligned}
    \left(\frac{1}{\tau_{k}g_{i,i}}I-J_{k}\right)Y_{k,i}=&
    ~\sum_{j=1}^{s}\frac{\bb_{i,j}}{\tau_{k}}Y_{k-1,j}+
    \sum_{j=1}^{s}a_{i,j}F(t_{k-1,j},\sum_{\ell=1}^{j}\bg_{j,\ell}Y_{k-1,\ell})\\
    &-J_{k}\sum_{j=1}^{s}\ba_{i,j}Y_{k-1,j}
    -\sum_{j=1}^{i-1}\frac{\bg_{i,j}}{\tau_{k}}Y_{k,j},\quad i=1,\dots,s.
  \end{aligned}\label{eq:auxpeerRos}
\end{align}
Again, replacing the Jacobian
$J_{k}$ by~(\ref{eq:frechet}), the modified Rosenbrock-type scheme, applied
to the DRE, reads
\begin{align}
  \begin{aligned}
    \tA_{k,i}^{T}Y_{k,i}&+Y_{k,i}\tA_{k,i}=-\tW_{k,i},\\
    \tW_{k,i}=&\sum_{j=1}^{s}\frac{\bb_{i,j}}{\tau_{k}}Y_{k-1,j}
    +\sum_{j=1}^{s}a_{i,j}\cR(t_{k-1,j},\sum_{\ell=1}^{j}\bg_{j,\ell}Y_{k-1,\ell})\\
    &-\sum_{j=1}^{s}\ba_{i,j}(\hA_{k}^{T}Y_{k-1,j}+Y_{k-1,j}\hA_{k})
    -\sum_{j=1}^{i-1}\frac{\bg_{i,j}}{\tau_{k}}Y_{k,j}
  \end{aligned}\label{eq:ALE_mRosPeer}
\end{align}
with $\hA_{k}$ from the original scheme and
$\tA_{k,i}=\hA_{k}-\frac{1}{2\tau_{k}g_{i,i}}I$.

Recall that the introduction of the auxiliary variables is capable of avoiding
the application of the Jacobian $J_{k}$ to the sum of current stage variables.
Still, the application remains for the sum of the previously determined peer
variables. Moreover, in contrast to the classical Rosenbrock methods, the
original solution approximations $X_{k,i}$ have to be reconstructed from the
auxiliary variables by (\ref{eq:origvar_peer}). That is, the reconstruction
doubles the online storage amount for storing the solution approximations
$X_{k,i}$ and the corresponding auxiliary variables $Y_{k,i},~i=1,\dots,s$,
during the runtime of the integration method.

Summarizing, the linearly implicit peer methods result in the solution of $s$
ALEs, just like the classical Rosenbrock methods~\cite{Men07,BenM13}, but
directly compute the sought for solutions, instead of additional stage
variables. Moreover, the additional stage variables from the classical
Rosenbrock methods have a low stage order and therefore the integration
procedures may suffer from order reduction. The computation of peer variables in
the Rosenbrock-type peer scheme, sharing the same accuracy and stability
properties, can overcome this well-known disadvantage~\cite{PodWS05} and again
allows us to use larger time steps compared to the classical Rosenbrock methods.
\section{Efficient Solution using Low-Rank Representations}
\label{sec:low-rank}
As mentioned in the introduction, for small-scale problems, the implicit and
Rosenbrock-type peer methods can directly be applied to the DRE, in general
resulting in dense solutions. Thus, the explicit computation of the solution is
not recommended for large-scale applications. Based on the observation that the
solution to the ALEs in the innermost iteration often is of low numerical
rank~\cite{AntSZ02,Gra04,Pen00a,TruV07}, the literature provides a number of
solution methods for large-scale ALEs based on low-rank versions of the
alternating directions implicit (ADI) iteration and Krylov subspace
methods. First developments considered a two-term LRCF of both ALE solution
philosophies. Most recent improvements can, e.g., be found
in~\cite{BenKS13a,BenKS13b,BenKS14,Kue16} and~\cite{JaiK94,StyS12,DruSZ14},
respectively. Three-term LRSIF based formulations of these solution strategies
have first been investigated for the more general case of Sylvester
equations~\cite{BenLT09}. The specific application to ALEs, is extensively
studied in~\cite{LanMS14,LanMS15,Lan17}. The latter factorization is of major
importance for the efficient solution of differential matrix equations. That is,
the LRSIF allows to avoid complex data and arithmetic, arising within the
classical low-rank two-term representation of the ALEs within the classical
implicit integration schemes of order $\geq 2$. Note that for the implicit and
Rosenbrock-type peer schemes complex data and arithmetic, in general, already
occur for order \(1\). The LRSIF has proven to show considerably better
performance with respect to computational timings and storage amount in most
applications. Note that there is some exceptions, see~\cite{LanMS15,Lan17} for
details. Still, for the numerical experiments in Section~\ref{sec:numerics}, the
algorithms used are restricted to the LRSIF based schemes. Moreover, we restrict
to implementations using the ADI iteration for the solution of the innermost
ALEs.

In order to exploit the low-rank phenomenon, a suitable low-rank representation
of the right hand sides of the ALEs~(\ref{eq:ALE_implPeer}) and
(\ref{eq:ALE_RosPeer})/(\ref{eq:ALE_mRosPeer}) within the implicit and linearly
implicit Rosenbrock-type peer schemes, respectively, has to be found. In what
follows, the LRSIF representations are presented. A detailed derivation of the
LRCF based strategy and an extension to generalized DREs, also for the LRSIF
approach, can be found in~\cite{Lan17}. For the remainder, we define the
mapping
\begin{align*}
  H:\R^{n\times n}\rightarrow \R^{2n\times 2n}, \quad H:I\mapsto H(I)=
  \begin{bmatrix}
    0&I\\
    I&0
  \end{bmatrix}.
\end{align*}

\subsection{Low-Rank Implicit Peer Scheme}
For the solution of the DRE~(\ref{eq:DRE}) by implicit peer schemes, the main
ingredient is to solve the algebraic Lyapunov equation
\begin{align}
  \begin{aligned}
    &{{}\hA_{k,i}^{(\ell)}}^TX_{k,i}^{(\ell)}+X_{k,i}^{(\ell)}\hA_{k,i}^{(\ell)}
    =-\tW_{k,i}-\tau_{k}g_{i,i}X_{k,i}^{(\ell-1)}S_{k,i}X_{k,i}^{(\ell-1)},\\
    &\tA_{k,i}=\tau_{k}g_{i,i}A_{k,i}-\frac{1}{2}I,\quad
    \tS_{k,i}=\tau_{k}g_{i,i}S_{k,i},\\
    &\hA_{k,i}^{(\ell)}=\tA_{k,i}-\tS_{k,i}X_{k,i}^{(\ell-1)}
    -\frac{1}{2}I,\\
    &\tW_{k,i}=\tau_{k}g_{i,i}W_{k,i}\!+\!\sum_{j=1}^{s}b_{i,j}X_{k-1,j}
    \!+\!\tau_{k}\!\sum_{j=1}^{s}\!a_{i,j}\cR(t_{k-1,j},X_{k-1,j})
    \!+\!\tau_{k}\!\sum_{j=1}^{i-1}\!g_{i,j}\cR(t_{k,j},X_{k,j})
  \end{aligned}\label{eq:RHS_implpeer}
\end{align}
within every Newton step $\ell$ at each time step $t_{k}\to t_{k+1}$.  Using
low-rank versions of the ADI method, this requires the right hand side to be
given in low-rank form as well. Provided $S_{k,i}$ and $W_{k,i}$ in the
DRE~(\ref{eq:DRE}) are given in the form
\begin{align*}
  S_{k,i}=B_{k,i}B_{k,i}^{T},\quad W_{k,i}=C_{k,i}^{T}C_{k,i}
\end{align*}
with $B_{k,i}\in\Rnm$ and $C_{k,i}\in\Rqn$, $m,q\ll n$, the right hand side of
the ALE can also be written in factored form.
Assume that the previous solution approximations $X_{k,j}$'s admit a
decomposition of the form $X_{k,j}=L_{k,j}D_{k,j}L_{k,j}^{T}$ with
$L_{k,j}\in\R^{n\times n_{L_{k,j}}},~D_{k,j}\in\R^{n_{L_{k,j}}\times
  n_{L_{k,j}}}$ such that the right hand side of~(\ref{eq:RHS_implpeer}) can be
written in the form $-G_{k,i}^{(\ell)}S_{k,i}^{(\ell)}{{}G_{k,i}^{(\ell)}}^{T}$.
In order to find such a symmetric indefinite decomposition of the entire right
hand side, we first define a factorization for the Riccati operator $\cR(.,.)$
in the form
\begin{align}
  \begin{aligned}
    \cR(t_{k,j},\!X_{k,j})&\!=\!C_{k,j}^{T}C_{k,j}\!+\!A_{k,j}^{T}X_{k,j}\!
    +\!\!X_{k,j}A_{k,j}
    \!-\!\!X_{k,j}B_{k,j}B_{k,j}^{T}X_{k,j}\!=\!\cT_{k,j}\cM_{k,j}\cT_{k,j}^{T},\\
    \cT_{k,j}&=
    \begin{bmatrix}
      C_{k,j}^{T},&A_{k,j}^{T}L_{k,j},&L_{k,j}
    \end{bmatrix}\in\R^{n\times (q+2n_{L_{k,j}})},\\
    \cM_{k,j}&=
    \begin{bmatrix}
      I_q&0&0\\
      0&0&D_{k,j}\\
      0&\quad D_{k,j}&\quad -D_{k,j}L_{k,j}^TB_{k,j}B_{k,j}^TL_{k,j}D_{k,j}\\
    \end{bmatrix}\!\in\R^{(q+2n_{L_{k,j}})\times (q+2n_{L_{k,j}})}.
  \end{aligned}\label{eq:RiccOp_peerLR}
\end{align}
For a more detailed derivation, we refer to~\cite{Lan17}. Then,
applying~(\ref{eq:RiccOp_peerLR}) to the right hand side
\begin{align*}
  -\tW_{k,i}&-\tau_{k}g_{i,i}X_{k,i}^{(\ell-1)}S_{k,i}X_{k,i}^{(\ell-1)},
\end{align*}
of the ALE from Equation~(\ref{eq:RHS_implpeer}), the decomposition
\(G_{k,i}^{(\ell)}S_{k,i}^{(\ell)}{{}G_{k,i}^{(\ell)}}^{T}\) is given by the
factors
\begin{align*}
  G_{k,i}^{(\ell)}=&
  \begin{bmatrix}
    C_{k,i}^{T},&L_{k-1,1},\dots,L_{k-1,s},&\cT_{k-1,1},\dots,\cT_{k-1,s},&
    \cT_{k,1},\dots,\cT_{k,i-1},&X_{k,i}^{(\ell-1)}B_{k,i}
  \end{bmatrix},\\
  S_{k,i}^{(\ell)}=&
    \diag\left(\tau_{k}g_{i,i}I_{q},~
      b_{i,1}D_{k-1,1},\dots,b_{i,s}D_{k-1,s},~
      \tau_{k}a_{i,1}\cM_{k-1,1},\dots,\tau_{k}a_{i,s}\cM_{k-1,s},\right.\\
       ~&\qquad\left.\tau_{k}g_{i,1}\cM_{k,1},\dots,\tau_{k}g_{i,i-1}\cM_{k,i-1},~
        \tau_{k}g_{i,i}I_{m}\right)
\end{align*}
can be formulated and the desired factor $G_{k,i}^{(\ell)}$ is of column size
\begin{align*}
  &~q+\sum_{j=1}^{s}n_{L_{k-1,j}}+\sum_{j=1}^{s}(q+2n_{L_{k-1,j}})
  +\sum_{j=1}^{i-1}(q+2n_{L_{k,j}})+n_{L_{k,i}^{(\ell-1)}}\\
  =&~(s+i)q+3\sum_{j=1}^{s}n_{L_{k-1,j}}+2\sum_{j=1}^{i-1}n_{L_{k,j}}+m.
\end{align*}
For autonomous systems with constant system matrices, the inner ALE becomes
\begin{align*}
  {{}\hA_{k,i}^{(\ell)}}^TX_{k,i}^{(\ell)}+X_{k,i}^{(\ell)}\hA_{k,i}^{(\ell)}
  =-\tW_{k,i}-\tau_{k}g_{i,i}X_{k,i}^{(\ell-1)}BB^{T}X_{k,i}^{(\ell-1)}=
  -G_{k,i}^{(\ell)}S_{k,i}^{(\ell)}{{}G_{k,i}^{(\ell)}}^{T},
\end{align*}
where $\hA_{k,i}$ and the right hand side factors
$G_{k,i}^{(\ell)},~S_{k,i}^{(\ell)}$ are given by
\begin{align*}
  \hA_{k,i}^{(\ell)}&=\!\!\tau_{k}g_{i,i}(A-BB^{T}X_{k,i}^{(\ell-1)})-\frac{1}{2}I\\
  G_{k,i}^{(\ell)}&=
  \begin{bmatrix}
    C^{T},&L_{k-1,1},\dots,L_{k-1,s},&\cT_{k-1,1},\dots,\cT_{k-1,s},&
    \cT_{k,1},\dots,\cT_{k,i-1},&X_{k,i}^{(\ell-1)}B
  \end{bmatrix},\\
  S_{k,i}^{(\ell)}=&
  \diag\left(\tau_{k}(\sum_{j=1}^{s}a_{i,j}+\sum_{j=1}^{i}g_{i,j})I_{q},~
                  b_{i,1}D_{k-1,1},\dots,b_{i,s}D_{k-1,s},\right.\\
  ~&~~\qquad\left.\tau_{k}a_{i,1}\cM_{k-1,1},\dots,\tau_{k}a_{i,s}\cM_{k-1,s},~
     \tau_{k}g_{i,1}\cM_{k,1},\dots,\tau_{k}g_{i,i-1}\cM_{k,i-1},~
     \tau_{k}g_{i,i}I_{m}\vphantom{\sum_{j=1}^{i}}\right)
\end{align*}
where the factors \(\cT_{k,j},\cM_{k,j}\) simplify to
\begin{align*}
  \cT_{k,j}&=
  \begin{bmatrix}
    A^{T}L_{k,j},&L_{k,j}
  \end{bmatrix}\in\R^{n\times 2n_{L_{k,j}}},\\
  \cM_{k,j}&=
  \begin{bmatrix}
    0&D_{k,j}\\
    D_{k,j}&\quad -D_{k,j}L_{k,j}^TB_{k,j}B_{k,j}^TL_{k,j}D_{k,j}\\
  \end{bmatrix}\in\R^{2n_{L_{k,j}}\times 2n_{L_{k,j}}}.
\end{align*}
Then, the right hand side factor $G_{k,i}^{(\ell)}$ is of column size
\begin{align*}
  q+3\sum_{j=1}^{s}n_{L_{k-1,j}}+2\sum_{j=1}^{i-1}n_{L_{k,j}}+m.
\end{align*}

\subsection{Low-Rank Rosenbrock-type Peer Scheme}
\subsubsection{Standard Rosenbrock-type Peer Representation}
\label{sec:LDL_RosPeer_DRE}
For the low-rank symmetric indefinite factorization based solution of a
non-autonomous DRE~(\ref{eq:DRE}), using the Rosenbrock-type peer method,
we consider the ALE
\begin{align}
  \begin{aligned}
    \tA_{k,i}^{T}X_{k,i}&+X_{k,i}\tA_{k,i}=-\tW_{k,i},\quad i=1,\dots,s,\\
    \tW_{k,i}=&~\sum_{j=1}^{s}b_{i,j}X_{k-1,j}
    +\tau_{k}\sum_{j=1}^{s}a_{i,j}\left(\cR(t_{k-1,j},X_{k-1,j})
      -(\hA_{k}^{T}X_{k-1,j}+X_{k-1,j}\hA_{k})\right),\\
    &+\tau_{k}\sum_{j=1}^{i-1}g_{i,j}(\hA_{k}^{T}X_{k,j}+X_{k,j}\hA_{k}),
  \end{aligned}\label{eq:stRosPeer_ALE}
\end{align}
where we have $\hA_{k}=A_{k}-B_{k}B_{k}^{T}X_{k},
~\tA_{k,i}=\tau_{k}g_{i,i}\hA_{k}-\frac{1}{2}I$. In contrast to
small-scale and dense computations, it is recommended to never explicitly form
the matrices $\hA_{k}$. Therefore, instead we use
\begin{align}
  \begin{aligned}
    \hA_{k}^{T}X_{k-1,j}+X_{k-1,j}\hA_{k}=&~A_{k}^{T}X_{k-1,j}+X_{k-1,j}A_{k}\\
    &-X_{k}B_{k}B_{k}^{T}X_{k-1,j}-X_{k-1,j}B_{k}B_{k}^{T}X_{k}.
  \end{aligned}\label{eq:lin_expand}
\end{align}
Using~(\ref{eq:lin_expand}) and further exploiting the structure of the Riccati
operators $\cR(t_{k-1,j},X_{k-1,j})$, the right hand side $\tW_{k,i}$ of the
standard Rosenbrock-type peer scheme~(\ref{eq:stRosPeer_ALE}) can be
reformulated in the form
\begin{align*}
  \begin{aligned}
    \tW_{k,i}=&~\tau_{k}\sum_{j=1}^{i-1}g_{i,j}\left(A_{k}^{T}X_{k,j}+X_{k,j}A_{k}
    -X_{k}B_{k}B_{k}^{T}X_{k,j}-X_{k,j}B_{k}B_{k}^{T}X_{k}\right)\\
    &+\sum_{j=1}^{s}\biggl(\tau_{k}a_{i,j}\bigl(C_{k-1,j}^{T}C_{k-1,j}
    -X_{k-1,j}B_{k-1,j}B_{k-1,j}^{T}X_{k-1,j}\biggr.\\
    &~\biggl.+X_{k}B_{k}B_{k}^{T}X_{k-1,j}+X_{k-1,j}B_{k}B_{k}^{T}X_{k}\bigr)
    +\check{A}_{k,i,j}^{T}X_{k-1,j}+X_{k-1,j}\check{A}_{k,i,j}\biggr),
  \end{aligned}
\end{align*}
where
$\check{A}_{k,i,j}=\tau_{k}a_{i,j}(A_{k-1,j}-A_{k})+\frac{b_{i,j}}{2}I$. The
matrix $\check{A}_{k,i,j}$ can efficiently be computed, since $A_{k-1,j}$
and $A_{k}$ are sparse matrices and so is $\check{A}_{k,i,j}$. Note that for $j=s$,
we have $A_{k-1,s}=A_{k},~B_{k-1,s}=B_{k}$ and $X_{k-1,s}=X_{k}$. Therefore
$\check{A}_{k,i,s}=\frac{b_{i,s}}{2}I$ and the right hand side at every stage
$i=1,\dots,s$ reduces to
\begin{align*}
  \begin{aligned}
    \tW_{k,i}=&~\tau_{k}\sum_{j=1}^{i-1}g_{i,j}\left(A_{k}^{T}X_{k,j}+X_{k,j}A_{k}
    -X_{k}B_{k}B_{k}^{T}X_{k,j}-X_{k,j}B_{k}B_{k}^{T}X_{k}\right)\\
    &+\sum_{j=1}^{s-1}\biggl(\tau_{k}a_{i,j}\bigl(C_{k-1,j}^{T}C_{k-1,j}
    -X_{k-1,j}B_{k-1,j}B_{k-1,j}^{T}X_{k-1,j}\biggr.\\
    &~\biggl.+X_{k}B_{k}B_{k}^{T}X_{k-1,j}+X_{k-1,j}B_{k}B_{k}^{T}X_{k}\bigr)
    +\check{A}_{k,i,j}^{T}X_{k-1,j}+X_{k-1,j}\check{A}_{k,i,j}\biggr)\\
    &+\tau_{k}a_{i,s}\left(C_{k}^{T}C_{k}+X_{k}B_{k}B_{k}^{T}X_{k}\right)
    +b_{i,s}X_{k}.
  \end{aligned}
\end{align*}
Also, we see that a considerable number of quadratic terms share the product
$X_{k}B_{k}$ or its transpose. Combining these expressions, we obtain the
formulation
\begin{align}
  \begin{aligned}
    \tW_{k,i}=&~\tau_{k}\sum_{j=1}^{i-1}g_{i,j}\left(A_{k}^{T}X_{k,j}+X_{k,j}A_{k}\right)
    +X_{k}B_{k}K_{k,i}^{T}+K_{k,i}B_{k}^{T}X_{k}\\
    &+\sum_{j=1}^{s-1}\biggl(\tau_{k}a_{i,j}\bigl(C_{k-1,j}^{T}C_{k-1,j}
    -X_{k-1,j}B_{k-1,j}B_{k-1,j}^{T}X_{k-1,j}\bigr)\\
    &+\check{A}_{k,i,j}^{T}X_{k-1,j}+X_{k-1,j}\check{A}_{k,i,j}\biggr)
    +\tau_{k}a_{i,s}C_{k}^{T}C_{k}+b_{i,s}X_{k},
  \end{aligned}\label{eq:rhs_RosPeer_DRE}
\end{align}
where
\begin{align*}
  K_{k,i}=\tau_{k}\left(\sum_{j=1}^{s-1}a_{i,j}X_{k-1,j}+\frac{a_{i,s}}{2}X_{k}
      -\sum_{j=1}^{i-1}g_{i,j}X_{k,j}\right)B_{k}
\end{align*}
collects all products, interacting with $X_{k}B_{k}$. Again, the previous
solution approximations $X_{k-1,j}=L_{k-1,j}D_{k-1,j}L_{k-1,j}^{T},j=1,\dots,s$,
$X_{k}=L_{k}D_{k}L_{k}^{T}$ with \(L_{k}=L_{k-1,s},D_{k}=D_{k-1,s}\) and
$X_{k,j}=L_{k,j}D_{k,j}L_{k,j}^{T},j=1,\dots,i-1$, are assumed to be given in
low-rank format. Then, defining the matrices
\begin{align*}
  \cT_{k,j}&=
  \begin{bmatrix}
    A_{k}^{T}L_{k,j},&L_{k,j}
  \end{bmatrix}\in\R^{n\times 2n_{L_{k,j}}},\ 
  \cM_{k,j}=\tau_{k}g_{i,j}H(D_{k,j})
  \in\R^{2n_{L_{k,j}} \times 2n_{L_{k,j}}},\\
  \check{\cT}_{k,i,j}&=
  \begin{bmatrix}
    C_{k-1,j}^{T},&\check{A}_{k,i,j}^{T}L_{k-1,j},&L_{k-1,j}
  \end{bmatrix},\\
  \check{\cM}_{k,i,j}&=
  \begin{bmatrix}
    \tau_{k}a_{i,j}I_{q}&0&0\\
    0&0&D_{k-1,j}\\
    0&\quad D_{k-1,j}&\quad
    -\tau_{k}a_{i,j}D_{k-1,j}L_{k-1,j}^{T}B_{k-1,j}B_{k-1,j}^{T}L_{k-1,j}D_{k-1,j}
  \end{bmatrix}
\end{align*}
with \(\check{\cT}_{k,i,j}\in\R^{n\times (q+2n_{L_{k-1,j}})}\),
\(\check{\cM}_{k,i,j}\in\R^{(q+2n_{L_{k-1,j}}) \times (q+2n_{L_{k-1,j}})}\), the
low-rank symmetric indefinite factorization
\(\tW_{k,i}=G_{k,i}S_{k,i}G_{k,i}^{T}\) of~(\ref{eq:rhs_RosPeer_DRE}) is given
by
\begin{align*}
  G_{k,i}=&
  \begin{bmatrix}
    \cT_{k,1},\dots,\cT_{k,i-1},&X_{k}B_{k},&K_{k,i},
    &\check{\cT}_{k,i,1},\dots,\check{\cT}_{k,i,s-1},&C_{k}^{T},&L_{k}
  \end{bmatrix},\\
  S_{k,i}=&\diag
    \left(\cM_{k,1},\dots,\cM_{k,i-1},~H(I_{m}),~
      \check{\cM}_{k,i,1},\dots,\check{\cM}_{k,i,s-1},~
      \tau_{k}a_{i,s}I_{q},~ b_{i,s}D_{k}\right)
\end{align*}
with $G_{k,i}$ being of column size
\begin{align}
  \begin{aligned}
    &~\sum_{j=1}^{i-1}2n_{L_{k,j}}+2m+\sum_{j=1}^{s-1}(q+2n_{L_{k-1,j}})+q+n_{L_{k}}\\
    =&~2\sum_{j=1}^{i-1}n_{L_{k,j}}+2\sum_{j=1}^{s-1}n_{L_{k-1,j}}+n_{L_{k}}+sq+2m.
  \end{aligned}\label{eq:clmsz_RosPeer}
\end{align}

In the autonomous case, we in particular have $A_{k-1,j}=A_{k}=A$. Hence,
$\check{A}_{k,i,j}=\frac{b_{i,j}}{2}I$, $i,j=1,\dots,s$, and together with the
modifications for $j=s,~X_{k-1,s}=X_{k}$, the right hand side $\tW_{k,i}$
in~(\ref{eq:rhs_RosPeer_DRE}) becomes
\begin{align*}
  \begin{aligned}
    \tW_{k,i}=&
    \tau_{k}\sum_{j=1}^{i-1}g_{i,j}\left(A_{k}^{T}X_{k,j}+X_{k,j}A_{k}\right)
    +X_{k}B_{k}K_{k,i}^{T}+K_{k,i}B_{k}^{T}X_{k}\\
    &+\sum_{j=1}^{s}\left(\tau_{k}a_{i,j}C^{T}C+b_{i,j}X_{k-1,j}\right)
    -\sum_{j=1}^{s-1}\tau_{k}a_{i,j}X_{k-1,j}BB^{T}X_{k-1,j}.
  \end{aligned}
\end{align*}
Then, similar to the non-autonomous scheme, for the simplified right hand side,
we have
\begin{align*}
  G_{k,i}=&
  \begin{bmatrix}
    \cT_{k,1},\dots,\cT_{k,i-1},&X_{k}B,&K_{k,i},
    &C^{T},&L_{k-1,1},\dots,L_{k-1,s-1},&L_{k}
  \end{bmatrix},\\
  S_{k,i}=&
            \diag\left(\cM_{k,1},\dots,\cM_{k,i-1},~ H(I_{m}),~
            \tau_{k}\sum_{j=1}^{s}a_{i,j}I_{q},~
            \tD_{k-1,1},\dots,\tD_{k-1,s-1},~ b_{i,s}D_{k}\right)
\end{align*}
where
\begin{align*}
  \cT_{k,j}&=
  \begin{bmatrix}
    A_{k}^{T}L_{k,j},&L_{k,j}
  \end{bmatrix}\in\R^{n\times 2n_{L_{k,j}}},\ 
  \cM_{k,j}=\tau_{k}g_{i,j}H(D_{k,j})
                        \in\R^{2n_{L_{k,j}} \times 2n_{L_{k,j}}},\\
  \tD_{k-1,j}&=b_{i,j}D_{k-1,j}-\tau_{k}a_{i,j}D_{k-1,j}L_{k-1,j}BB^{T}L_{k-1,j}D_{k-1,j}.
\end{align*}
Here, the column size of the factor $G_{k,i}$ is
\begin{align}
  \sum_{j=1}^{i-1}2n_{L_{k,j}}+2m+q+\sum_{j=1}^{s}n_{L_{k-1,j}}
  =2\sum_{j=1}^{i-1}n_{L_{k,j}}+\sum_{j=1}^{s}n_{L_{k-1,j}}+q+2m.
  \label{eq:clmsz_RosPeer_aut}
\end{align}
\subsubsection{Modified Rosenbrock-type Peer Representation}
\label{sec:LDL_mRosPeer_DRE}
Now, for the modified Rosenbrock-type peer formulation applied to the
non-autonomous DRE, we consider the ALE
\begin{align*}
  \begin{aligned}
    \tA_{k,i}^{T}Y_{k,i}&+Y_{k,i}\tA_{k,i}=-\tW_{k,i},\quad i=1,\dots,s,\\
    \tW_{k,i}=&\sum_{j=1}^{s}\frac{\bb_{i,j}}{\tau_{k}}Y_{k-1,j}
    +\sum_{j=1}^{s}a_{i,j}\cR(t_{k-1,j},\sum_{\ell=1}^{j}\bg_{j,\ell}Y_{k-1,\ell})\\
    &-\sum_{j=1}^{s}\ba_{i,j}(\hA_{k}^{T}Y_{k-1,j}+Y_{k-1,j}\hA_{k})
    -\sum_{j=1}^{i-1}\frac{\bg_{i,j}}{\tau_{k}}Y_{k,j},\\
    \tA_{k,i}=&~\hA_{k}-\frac{1}{2\tau_{k}g_{i,i}}I,\quad
    \hA_{k}=A_{k}-B_{k}B_{k}^{T}X_{k}I.
  \end{aligned}
\end{align*}
Note that the matrix $\hA_{k}$ is given in terms of $X_{k}$. This is due to the
fact that $\hA_{k}$ originates from the Jacobian~(\ref{eq:frechet}) that, as in
the original scheme, is given as the Fr\'{e}chet derivative of
$\cR(t_{k},X_{k})=\cR(t_{k-1,s},X_{k-1,s})=
\cR(t_{k-1,s},\sum_{\ell=1}^{s}\bg_{j,\ell}Y_{k-1,\ell})$.  Thus, instead of
explicitly forming $\hA_{k}$, again relation~(\ref{eq:lin_expand}) is utilized.

For the sake of simplicity the original variables $X_{k}$ within $\hA_{k}$, as
well as in the Riccati operators $\cR(t_{k-1,j},X_{k-1,j})$, are kept throughout
the computations. As previously mentioned in Section~\ref{sec:RosPeer}, we have
to reconstruct the solution from the auxiliary variables anyway. Thus, using
both sets of variables does not require additional computations. In order to
give a more detailed motivation for mixing up the original and auxiliary scheme,
the following considerations are stated.

From the relation of the original and auxiliary variables, given
in~(\ref{eq:origvar_peer}), we have
\begin{align*}
X_{k-1,j}=\sum_{\ell=1}^{j}\bg_{j,\ell}Y_{k-1,\ell}.
\end{align*}
Further, defining the decomposition
$Y_{k-1,\ell}=\hL_{k-1,\ell}\hD_{k-1,\ell}\hL_{k-1,\ell}^{T}$, $\ell=1,\dots,j$
with $\hL_{k-1,\ell}\in\R^{n\times n_{\hL_{k-1,\ell}}}$,
\(\hD_{k-1,\ell}\in\R^{n_{\hL_{k-1,\ell}}\times n_{\hL_{k-1,\ell}}}\), the
original solution approximation admits a factorization
$X_{k-1,j}= L_{k-1,j}D_{k-1,j}L_{k-1,j}^{T}$, $j=1,\dots,s$, based on the
factors
\begin{align*}
  L_{k-1,j}&=
  \begin{bmatrix}
    \hL_{k-1,1},\dots,\hL_{k-1,j}
  \end{bmatrix},\quad
  D_{k-1,j}=\diag{\bg_{j,1}\hD_{k-1,1},\dots,\bg_{j,j}\hD_{k-1,j}}.
\end{align*}
The factors
$L_{k-1,j}\in\R^{n\times n_{L_{k-1,j}}},D_{k-1,j}\in\R^{n_{L_{k-1,j}}\times
  n_{L_{k-1,j}}}$ are given as a block concatenation of the solution factors of
the auxiliary variables $Y_{k-1,\ell}$, $\ell=1,\dots,j$. That is, the column
size $n_{L_{k-1,j}}=\sum_{\ell=1}^{j}n_{\hL_{k-1,\ell}}$ may dramatically grow
with respect to the number of stages and time steps. Still, the numerical rank
of the original solution is assumed to be ``small''. Thus, using column
compression techniques, see~\cite[Section 6.3]{Lan17}, being a tacit requirement
for large-scale problems anyway, the column size of $L_{k-1,j}$ is presumably
``small'' as well. To be more precise, the factors $L_{k-1,j}$ and $\hL_{k-1,j}$
are expected to be of compatible size. Consequently, one can make use of both
representations at the one place or another without messing up the formulations
with respect to both, the notational and computational complexity.

However, expanding $\hA_{k}$ and combining the linear parts with respect
to $Y_{k-1,j}$, the right hand side reads
\begin{align}
  \begin{aligned}
    \tW_{k,i}\!=&\!
    -\!\!\sum_{j=1}^{s}\biggl(\check{A}_{k,i,j}^{T}Y_{k-1,j}+Y_{k-1,j}\check{A}_{k,i,j}
    \!-\!\ba_{i,j}\left(X_{k}B_{k}B_{k}^{T}Y_{k-1,j}
    \!+\!Y_{k-1,j}B_{k}B_{k}^{T}X_{k}\right)\biggr)\\
    &+\sum_{j=1}^{s}a_{i,j}\cR(t_{k-1,j},X_{k-1,j})
    -\sum_{j=1}^{i-1}\frac{\bg_{i,j}}{\tau_{k}}Y_{k,j}
  \end{aligned}\label{eq:rhs_mRosPeer_DRE}
\end{align}
with $\check{A}_{k,i,j}=\ba_{i,j}A_{k}-\frac{\bb_{i,j}}{2\tau_{k}}I$.  Then,
separating $\cR(t_{k-1,s},X_{k-1,s})=\cR(t_{k},X_{k})$ and again combining the
quadratic terms including the products $X_{k}B_{k}$, we end up with
the formulation
\begin{align*}
  \begin{aligned}
    \tW_{k,i}=&-\sum_{j=1}^{s}
    \biggl(\check{A}_{k,i,j}^{T}Y_{k-1,j}+Y_{k-1,j}\check{A}_{k,i,j}\biggr)
    +X_{k}B_{k}K_{k,i}^{T}+K_{k,i}B_{k}^{T}X_{k}\\
    &+\sum_{j=1}^{s-1}a_{i,j}\cR(t_{k-1,j},X_{k-1,j})
    +a_{i,s}\left(C_{k}^{T}C_{k}+A_{k}^{T}X_{k}+X_{k}A_{k}\right)
    -\sum_{j=1}^{i-1}\frac{\bg_{i,j}}{\tau_{k}}Y_{k,j},\\
    K_{k,i}=&~\left(
      \sum_{j=1}^{s}\ba_{i,j}Y_{k-1,j}-\frac{a_{i,s}}{2}X_{k}
    \right)B_{k}.
  \end{aligned}
\end{align*}

Then, the associated symmetric indefinite formulation is given by the factors
\begin{align*}
  G_{k,i}=&~\bigl[
  \check{\cT}_{k-1,i,1},\dots,\check{\cT}_{k-1,i,s},
  ~X_{k}B_{k},~K_{k,i},
  ~\cT_{k-1,1},\dots,\cT_{k-1,s-1},\\
  &~~~C_{k}^{T},A_{k}^{T}L_{k},~L_{k},\sqrt{a_{i,s}}
    \hL_{k,1},\dots,\hL_{k,i-1},
  \bigr],\\
  S_{k,i}=&
            \diag\left(\vphantom{\frac{\bg_{i,i-1}}{\tau_{k}}}
            -\check{\cM}_{k-1,i,1},\dots,-\check{\cM}_{k-1,i,s},~ 
            H(I_{m}),~ a_{i,1}\cM_{k-1,1},\dots,a_{i,s-1}\cM_{k-1,s-1},\right.\\
          ~&\qquad~ \left. a_{i,s}I_{q},~ a_{i,s}H(D_{k}),~
            -\frac{\bg_{i,1}}{\tau_{k}}\hD_{k,1},\dots,
            -\frac{\bg_{i,i-1}}{\tau_{k}}\hD_{k,i-1}
            \right)
\end{align*}
with
\begin{align*}
  \check{\cT}_{k-1,i,j}&=
  \begin{bmatrix}
    \check{A}_{k,i,j}^{T}\hL_{k-1,j},&\hL_{k-1,j}
  \end{bmatrix}\!\in\R^{n\times 2n_{\hL_{k-1,j}}},
  \check{\cM}_{k-1,i,j}=H(\hD_{k-1,j})\!
  \in\R^{2n_{\hL_{k-1,j}}\times 2n_{\hL_{k-1,j}}},\\
  \cT_{k-1,j}&=
  \begin{bmatrix}
    C_{k-1,j}^{T},&A_{k-1,j}^{T}L_{k-1,j},&L_{k-1,j}
  \end{bmatrix}\in\R^{n\times (q+2n_{L_{k-1,j}})},\\
  \cM_{k-1,j}&=
  \begin{bmatrix}
    I_{q}&\!\! 0&\!\! 0\\
    0&\!\! 0&\!\!D_{k-1,j}\\
    0 &\!\! D_{k-1,j}&\!\!
    -D_{k-1,j}L_{k-1,j}^{T}B_{k-1,j}B_{k-1,j}^{T}L_{k-1,j}D_{k-1,j}
  \end{bmatrix}\!\!\!\in\!\R^{(q+2n_{L_{k-1,j}})\times(q+2n_{L_{k-1,j}})},
\end{align*}
defining the factorization of the Lyapunov-type expression and the Riccati
operators, respectively. The resulting column size of $G_{k,i}$ is then given by
\begin{align}
  \begin{aligned}
    &~\sum_{j=1}^{s}2n_{\hL_{k-1,j}}+2m+\sum_{j=1}^{s-1}(q+2n_{L_{k-1,j}})
    +q+2n_{L_{k}}+\sum_{j=1}^{i-1}n_{\hL_{k,j}}\\
    =&~\sum_{j=1}^{i-1}n_{\hL_{k,j}}
    +2\sum_{j=1}^{s}(n_{\hL_{k-1,j}}+n_{L_{k-1,j}})+sq+2m.
  \end{aligned}\label{eq:clmsz_mRosPeer}
\end{align}
Note that the use of both, the auxiliary variables in the linear parts and the
original variables within the Fr\'{e}chet derivative and the Riccati operator,
does not allow us to completely combine these parts, as we have seen for the
condensed form~(\ref{eq:rhs_RosPeer_DRE}) of the original Rosenbrock-type peer
scheme. Therefore, assume the associated low-rank factors $L_{k-1,j},D_{k-1,j}$
and $\hL_{k-1,j},\hD_{k-1,j}$ of $X_{k-1,j}$ and $Y_{k-1,j}$, respectively, to
be of comparable column sizes $n_{L_{k-1,j}}$ and $n_{\hL_{k-1,j}}$. Then,
comparing (\ref{eq:clmsz_RosPeer}) and~(\ref{eq:clmsz_mRosPeer}), the modified
scheme results in a larger overall number of columns in the right hand side
factorization, although avoiding the application of the Jacobian to the current
solutions $Y_{k,j},~j=1,\dots,i$, saves $2\sum_{j=1}^{i-1}n_{\hL_{k,j}}$ columns
in the first place. That is, for large-scale non-autonomous DREs, the standard
version of the Rosenbrock-type peer schemes seems to be preferable.

Still, a more beneficial situation can be found for autonomous DREs. Here,
additional modifications, based on the time-invariant nature of the system
matrices, allow to further reduce the complexity of the ALEs to be solved. In
that case, the associated ALEs are of the form
\begin{align}
  \begin{aligned}
    \tA_{k,i}^{T}Y_{k,i}&+Y_{k,i}\tA_{k,i}=-\tW_{k,i},\quad i=1,\dots,s,\\
    \tW_{k,i}=&\sum_{j=1}^{s}\frac{\bb_{i,j}}{\tau_{k}}Y_{k-1,j}
    +\sum_{j=1}^{s}a_{i,j}\cR(\sum_{\ell=1}^{j}\bg_{j,\ell}Y_{k-1,\ell})\\
    &-\sum_{j=1}^{s}\ba_{i,j}(\hA_{k}^{T}Y_{k-1,j}+Y_{k-1,j}\hA_{k})
    -\sum_{j=1}^{i-1}\frac{\bg_{i,j}}{\tau_{k}}Y_{k,j},\\
    \tA_{k,i}=&~\hA_{k}-\frac{1}{2\tau_{k}g_{i,i}},\quad
    \hA_{k}=A-BB^{T}X_{k}.
  \end{aligned}\label{eq:mRosPeer_autDRE}
\end{align}
We start the investigations at
$\cR(\sum_{\ell=1}^{j}\bg_{j,\ell}Y_{k-1,\ell})$. For that, first consider the
sum of Riccati operators
\begin{align*}
  \sum_{j=1}^{s}a_{i,j}\cR(X_{k-1,j})
  =\sum_{j=1}^{s}a_{i,j}\left(C^{T}C+A^{T}X_{k-1,j}+X_{k-1,j}A
    -X_{k-1,j}BB^{T}X_{k-1,j}\right).
\end{align*}
Further, recall the definitions $\bX_{k}=(X_{k,i})_{i=1}^{s}$ and
$\bY_{k}=(Y_{k,i})_{i=1}^{s}$. Then, from $A_{k}=A$ being constant and motivated
by~(\ref{eq:modsumPeer}), for the linear part, we find
\begin{align*}
  \begin{aligned}
    &\sum_{j=1}^{s}a_{i,j}A^{T}X_{k-1,j}+\sum_{j=1}^{s}a_{i,j}X_{k-1,j}A,~i=1,\dots,s\\
    \Leftrightarrow\quad& ((a_{i,j})\otimes A^{T})\bX_{k-1}+((a_{i,j})\otimes I
    )\bX_{k-1}A.
  \end{aligned}
\end{align*}
Moreover, from the definition~(\ref{eq:origvar_peer}) of $\bX_{k}$ in terms of
the auxiliary variables $\bY_{k}$ the following reformulation holds:
\begin{align*}
  \begin{aligned}
    &~((a_{i,j})\otimes A^{T})\bX_{k-1}\!+\!((a_{i,j})\otimes I)\bX_{k-1}A\\
    =&~((a_{i,j})\otimes A^{T})(G^{ -1}\otimes I)\bY_{k-1}
    +((a_{i,j})\otimes I)(G^{-1}\otimes I)\bY_{k-1}A\\
    =&~((a_{i,j})G^{-1}\otimes A^{T})\bY_{k-1}
    +((a_{i,j})G^{-1}\otimes I)\bY_{k-1}A\\
    =&~((\ba_{i,j})\otimes A^{T})\bY_{k-1}+((\ba_{i,j})\otimes I
    )\bY_{k-1}A.
  \end{aligned}
\end{align*}
with $(\ba_{i,j})=(a_{i,j})G^{-1}$ from~(\ref{auxcoeff_peer}). Then, together with
\begin{align*}
  \begin{aligned}
    &\quad((\ba_{i,j})\otimes A^{T})\bY_{k-1}+((\ba_{i,j})\otimes I)\bY_{k-1}A\\
    \Leftrightarrow&\sum_{j=1}^{s}\ba_{i,j}A^{T}Y_{k-1,j}+
    \sum_{j=1}^{s}\ba_{i,j}Y_{k-1,j}A,~i=1,\dots,s,
  \end{aligned}
\end{align*}
the sum of Riccati operators $\cR(X_{k-1,j})$ can be written in the mixed form
\begin{align*}
  \sum_{j=1}^{s}a_{i,j}\cR(X_{k-1,j})
  &=\sum_{j=1}^{s}a_{i,j}\left(C^{T}C+A^{T}X_{k-1,j}+X_{k-1,j}A
    -X_{k-1,j}BB^{T}X_{k-1,j}\right)\\
  &=\sum_{j=1}^{s}a_{i,j}\left(C^{T}C-X_{k-1,j}BB^{T}X_{k-1,j}\right)
  +\sum_{j=1}^{s}\ba_{i,j}\left(A^{T}Y_{k-1,j}+Y_{k-1,j}A\right).
\end{align*}
Note that in this formulation only the quadratic term of the Riccati operator
uses the original variables and analogously to the right hand side $\tW_{k,i}$
in~(\ref{eq:rhs_mRosPeer_DRE}), for an autonomous DRE, we obtain
\begin{align*}
  \begin{aligned}
    \tW_{k,i}=&
    -\sum_{j=1}^{s}\biggl(\check{A}_{k,i,j}^{T}Y_{k-1,j}+Y_{k-1,j}\check{A}_{k,i,j}
    -\ba_{i,j}\left(X_{k}BB^{T}Y_{k-1,j}
    +Y_{k-1,j}BB^{T}X_{k}\right)\biggr)\\
    &+\sum_{j=1}^{s}a_{i,j}\left(C^{T}C-X_{k-1,j}BB^{T}X_{k-1,j}\right)
    +\sum_{j=1}^{s}\ba_{i,j}\left(A^{T}Y_{k-1,j}+Y_{k-1,j}A\right)\\
    &-\sum_{j=1}^{i-1}\frac{\bg_{i,j}}{\tau_{k}}Y_{k,j}
  \end{aligned}
\end{align*}
with $\check{A}_{k,i,j}=\ba_{i,j}A-\frac{\bb_{i,j}}{2\tau_{k}}I$. Now, combining
the expressions that are linear in $Y_{k-1,j}$, as well as the quadratic terms
containing $X_{k}B$ and again paying particular attention to $j=s$ with
$X_{k-1,s}=X_{k}$, $Y_{k-1,s}=Y_{k}$, the right hand side reads
\begin{align}
  \begin{aligned}
    \tW_{k,i}=&
    ~\sum_{j=1}^{s}a_{i,j}C^{T}C-\sum_{j=1}^{s-1}a_{i,j}X_{k-1,j}BB^{T}X_{k-1,j}
    +X_{k}BK_{k,i}^{T}+K_{k,i}B^{T}X_{k}\\
    &+\sum_{j=1}^{s}\frac{\bb_{i,j}}{\tau_{k}}Y_{k-1,j}
    -\sum_{j=1}^{i-1}\frac{\bg_{i,j}}{\tau_{k}}Y_{k,j},\\
    K_{k,i}=&~\left(
      \sum_{j=1}^{s}\ba_{i,j}Y_{k-1,j}-\frac{a_{i,s}}{2}X_{k}
    \right)B.
  \end{aligned}\label{eq:rhs_mRosPeer_autDRE}
\end{align}

For the autonomous case and the associated ALE~(\ref{eq:mRosPeer_autDRE}) and
its condensed right hand side~(\ref{eq:rhs_mRosPeer_autDRE}), we find the
factors
\begin{align*}
  G_{k,i}=&~\bigl[
    C^{T},~X_{k-1,1}B,\dots,X_{k-1,s-1}B,~X_{k}B,~K_{k,i},\\
    &~~~\hL_{k-1,1},\dots,\hL_{k-1,s},~\hL_{k,1},\dots,\hL_{k,i-1}
  \bigr],\\
  S_{k,i}=&
            \diag\left(
            \sum_{j=1}^{s}a_{i,j}I_{q},~ -a_{i,1}I_{m},\dots,-a_{i,s-1}I_{m},~
            H(I_{m}),\right.\\
          ~&\qquad~ \left.
            \frac{\bb_{i,1}}{\tau_{k}}\hD_{k-1,1},\dots,
            \frac{\bb_{i,s}}{\tau_{k}}\hD_{k-1,s},~
            -\frac{\bg_{i,1}}{\tau_{k}}\hD_{k,1},\dots,
            -\frac{-\frac{1}{\tau_{k}}}{\tau_{k}}\hD_{k,i-1}
            \right)
\end{align*}
where $G_{k,i}$ is of column size
\begin{align}
  q+\sum_{j=1}^{s-1}m+2m+\sum_{j=1}^{s}n_{\hL_{k-1,j}}+\sum_{j=1}^{i-1}n_{\hL_{k,j}}
  =\sum_{j=1}^{i-1}n_{\hL_{k,j}}+\sum_{j=1}^{s}n_{\hL_{k-1,j}}+q+(s+1)m.
  \label{eq:clmsz_mRosPeer_aut}
\end{align}
Again, assume that the column sizes of the solution factors \(L_{k,j}\) and
\(\hL_{k,j}\) of the original Rosenbrock-type peer and its modified version,
respectively, are compatible. Then, from (\ref{eq:clmsz_RosPeer_aut}) and
(\ref{eq:clmsz_mRosPeer_aut}) it can be seen that the modified version can save
a number of system solves within the ALE solver, as long as \((s-1)m\) does not
exceed \(\sum_{j=1}^{i-1}n_{L_{k,j}}\) from the original scheme. This will most
likely be true for a small number \(m\), i.e, a low numerical rank of \(S(t)\)
in the DRE~\eqref{eq:DRE}. Considering control problems, \(m\) represents the
number of inputs to the system to be controlled and thus will be rather small
for numerous examples.
\section{Numerical Experiments}\label{sec:numerics}
The following computations have been executed on a 64bit CentOS 5.5 system with
two {Intel\textsuperscript{\textregistered}}\
{Xeon\textsuperscript{\textregistered}}\ X5650@2.67 GHz with a total of 12 cores
and 48GB main memory, being one computing node of the linux cluster {\rmfamily
  otto}\footnote{\url{http://www.mpi-magdeburg.mpg.de/1012477/otto}} at the
{\itshape Max Planck Institute for Dynamics of Complex Technical Systems} in
Magdeburg. The numerical algorithms have been implemented and tested in
{\matlab}~version 8.0.0.783 (R2012b).

For the numerical experiments, we consider the implicit peer
method~(\ref{eq:peer2}) and both versions of the Rosenbrock-type
schemes~(\ref{eq:peerRos}),~(\ref{eq:auxpeerRos}) up to order \(4\). For a
comparison of the computational times and relative errors with respect to a
reference solution, the several peer schemes are also compared to the BDF
methods of order 1 to 4~\cite{BenM04,LanMS15,Lan17}, Rosenbrock methods of orders
\(1,2\)~\cite{BenM13}, \(4\)~\cite{Sha82}, and the midpoint and trapezoidal
rules~\cite{Die92}. The relative errors are given in the Frobenius norm
$\|.\|_F$. An overview of the corresponding low-rank formulations, except for
the Rosenbrock method of order 4, can be found in~\cite{LanMS15,Lan17}. For the
latter, no low-rank representation has been published so far. The additional
initial values for multi-step and the peer integrators of order\(\geq 2\), the
one-step Rosenbrock methods of appropriate order are chosen. In what remains,
the abbreviations, given in Table~\ref{tab:acronyms}, are used to identify the
several integration schemes. For the integration methods, using Newton's method
to solve the arising AREs, a tolerance of \(1e\)-10 and a maximum number of 15
Newton steps are chosen. The ADI iteration, used in the innermost loop of all
schemes, is terminated at a tolerance of \(n\varepsilon\) or at a maximum of 100
ADI steps. Here again, \(n\) is the system dimension and \(\varepsilon\) denotes
the machine precision.
\begin{table}[t]
  \centering
  \caption{Acronyms of the time integration methods (\(s=1,\dots,4\)).}
  \label{tab:acronyms}
  {\rowcolors{2}{gray!25}{white}
    \begin{tabular}{cc}
      \rowcolor{gray!60}
      {Time integration method}&{Acronym}\\
      BDF of order \(s\)& BDF\((s)\)\\
      Rosenbrock of order \(s\) & Ros\((s)\)\\
      Midpoint rule & Mid\\
      Trapezoidal rule & Trap\\
      Implcit peer of order \(s\)& Peer\((s)\)\\
      Rosenbrock-type peer of order \(s\)& RosPeer\((s)\)\\
      Modified RosPeer\((s)\) & mRosPeer\((s)\)\\
    \end{tabular}
  }
\end{table}
\begin{table}[t]
  \centering
  \caption{$2$-stage implicit peer method of order $2$.}
  \label{tab:peer2}
  {\rowcolors{1}{white}{gray!25}
  \begin{tabular}{|lrlr|}\firsthline
    $c_{1}:$&$0.4831632475943920$ & $c_{2}:$&$1.0000000000000000$\\
    $b_{1,1}:$&$-0.3045407685048590$ & $b_{1,2}:$&$1.3045407685048591$\\
    $b_{2,1}:$&$-0.3045407685048590$ & $b_{2,2}:$&$1.3045407685048591$\\
    $g_{1,1}:$&$0.2584183762028040$ & $g_{1,2}:$&$0.0000000000000000$\\
    $g_{2,1}:$&$0.4376001712448750$ & $g_{2,2}:$&$0.2584183762028040$\\\hline
  \end{tabular}
  }
\end{table}

\paragraph{Implicit Peer Coefficients}
\label{sec:impl-peer-coeff}
The \(1\)-stage implicit peer scheme is given by the coefficients $c_{1}=1$,
$b_{1,1}=1$ and $g_{1,1}=1$. The coefficients of the $2$-stage implicit peer
method, given in Table~\ref{tab:peer2}, were provided by the group of
Prof. R. Weiner at the Martin-Luther-Universität Halle and cannot, to the best
of the authors' knowledge, be found in
any publication so far. The coefficients for the \(3\)- and \(4\)-stage peer
schemes are provided by methods 3a and 4b in~\cite{SolW17}. 

\paragraph{Rosenbrock-type Peer Coefficients}
\label{sec:rosenbrock-type-peer}
The \(1\)-stage Rosenbrock-type peer method is given by the coefficients
$c_{1}=1$, $a_{1,1}=1$, $b_{1,1}=1$ and $g_{1,1}=1$. The coefficients for the
Rosenbrock-type peer schemes used here, can be computed following the
instructions in~\cite[Section~3]{PodWS05}.

\subsection{Steel Profile}
\label{sec:steel-profile}
\begin{figure}[t]
  \centering
  \begin{subfigure}{.325\linewidth}
    \includegraphics{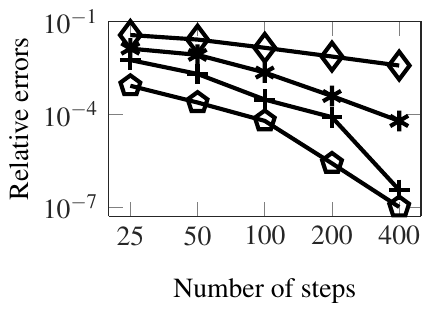}%
    \caption{Peer(1-4)}\label{fig:acc_Peer}
  \end{subfigure}\hskip-.25em%
  \begin{subfigure}{.325\linewidth}\hskip-.5em%
    \includegraphics{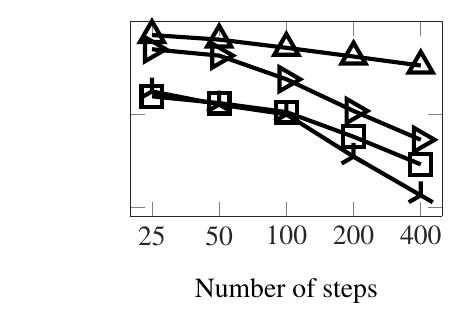}%
    \caption{RosPeer(1-4)}\label{fig:acc_RosPeer}
  \end{subfigure}\hskip-.25em%
  \begin{subfigure}{.325\linewidth}\hskip-.75em%
    \includegraphics{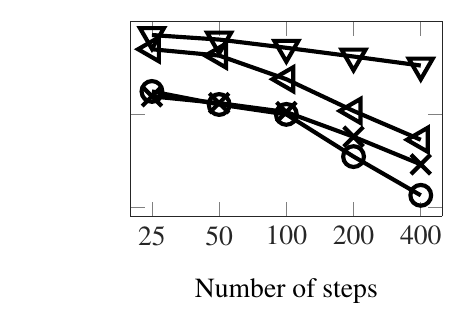}
    \caption{mRosPeer(1-4)}
    \label{fig:acc_modRosPeer}
  \end{subfigure}\vskip.25em%
  \begin{subfigure}{.325\linewidth}
    \includegraphics{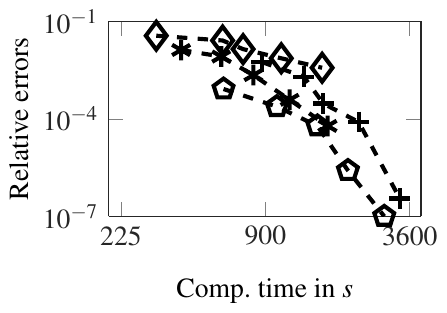}
    \caption{Peer(1-4)}\label{fig:eff_Peer}
  \end{subfigure}\hskip-.5em%
  \begin{subfigure}{.325\linewidth}\hskip-.5em
    \includegraphics{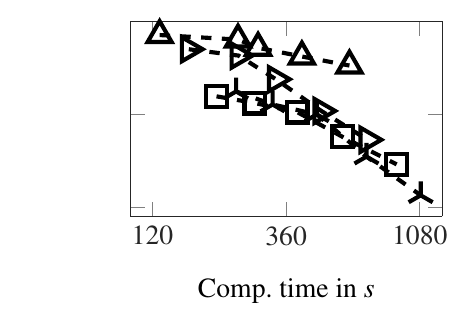}
    \caption{RosPeer(1-4)}\label{fig:eff_RosPeer}
  \end{subfigure}\hskip-.25em%
  \begin{subfigure}{.325\linewidth}\hskip-.75em
    \includegraphics{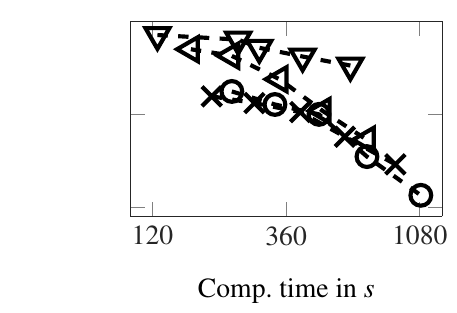}
    \caption{mRosPeer(1-4)}\label{fig:eff_modRosPeer}
  \end{subfigure}%
  \medskip
  \begin{center}\vskip-1em
    \includegraphics[scale=1.25]{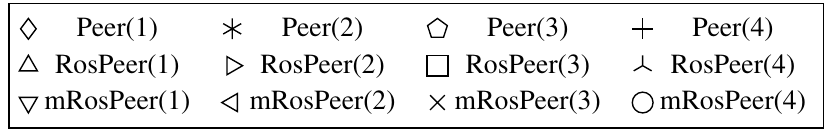}
  \end{center}\vskip-1em
  \caption{Steel profile: Accuracy and efficiency plots}
  \label{fig:acc_eff}
\end{figure}

{\rowcolors{3}{white}{gray!25}
  \begin{table}[t]
    \centering
    \caption{Steel profile: Computational timings and relative errors with
      respect to the reference solution for $\tau=0.1125\ s$, 400
      steps.}\label{tab:rail_timings}
    \begin{tabular}{lrc}
      \rowcolor{gray!50}
      Method&Time in $s$&Rel. Frobenius err.\\
      BDF(1)&1\,627.76&3.75e-03\\
      BDF(2)&1\,347.55&3.20e-04\\
      BDF(3)&1\,228.07&1.29e-04\\
      BDF(4)&1\,179.00&4.58e-05\\
      Ros1&806.00&3.75e-03\\
      Ros2&1\,028.04&1.24e-03\\
      Ros4&1\,001.05&1.30e-06\\
      Mid&1\,239.78&1.33e-04\\
      Trap&1\,202.90&1.32e-04\\
      Peer(1)&1\,551.03&3.75e-03\\
      Peer(2)&1\,635.30&6.09e-05\\
      Peer(3)&2\,815.84&1.01e-07\\
      Peer(4)&3\,268.56&3.57e-07\\
      RosPeer(1)&605.64&3.75e-03\\
      RosPeer(2)&702.46&1.50e-05\\
      RosPeer(3)&892.35&2.41e-06\\
      RosPeer(4)&1\,087.55&2.41e-07\\
      mRosPeer(1)&610.31&3.75e-03\\
      mRosPeer(2)&698.74&1.50e-05\\
      mRosPeer(3)&883.33&2.41e-06\\
      mRosPeer(4)&1\,088.86&2.41e-07
    \end{tabular}
  \end{table}
}
As a first example, we consider a semi-discretized heat transfer problem for
optimal cooling of steel profiles~\cite{BenS05b,morwiki_steel}. This example is
a mutli-input multi-output (MIMO) system with $m=7$ inputs and $q=6$
outputs. The solution to the DRE is computed on the simulation time interval
$[0,\;4\,500]\ s$ with the step sizes $\tau\in\{180, 90, 45, 25.5, 12.75\}\ s$
and \(\{25, 50, 100, 200, 400\}\) steps, respectively. Note that the actual time
line is implicitly scaled by $1e2$ within the model such that a real time of
$[0,45]\ s$ with corresponding step sizes is investigated. To ensure the
computability of a reference solution in appropriate time, the smallest
available discretization level with $n=371$ is chosen. The reference is computed
by the small-scale dense version of the fourth-order Rosenbrock (Ros4)
method. In particular, the {\itshape Parareal} based implementation with \(450\)
coarse and additionally \(1000\) fine steps at each of those intervals,
considered in~\cite{KoeLS16}, has been used.

Figures~\ref{fig:acc_eff}(\subref{fig:acc_Peer})-(\subref{fig:acc_modRosPeer})
show the accuracy plots for the implicit peer methods, the RosPeer schemes and
the modified RosPeer integrators, respectively. It can be observed that, for
this example, the convergence orders are reached asymptotically. Further, note
that the Peer(3) scheme outperforms its Peer(4) successor. This is due to the
superconvergence of the Peer(3) method (see~\cite[Section 4, Method 3a]{SolW17})
and the fact that, for this example, the convergence order 4 of the Peer(4)
scheme has just been reached for the last step size refinement. It can further
be observed that the implicit peer and the Rosenbrock-type schemes of
corresponding order achieve a comparable accuracy. This is not too surprising
considering an autonomous problem. The efficiency plots are presented in
Figures~(\subref{fig:eff_Peer})-(\subref{fig:eff_modRosPeer}). In
Table~\ref{tab:rail_timings}, the LRSIF computation times and the relative
errors with respect to the reference solution are given. Here, it becomes clear
that the peer methods of order \(s\geq 2\) show a significantly better
performance compared to the other implicit time integrators of similar order
with respect to the accuracy. Solely comparing the computational times, the
Rosenbrock-type peer scheme of first-order shows best performance. Taking the
efficiency into account, i.e., studying the required computational time versus
the achieved error level, see also
Figures~\ref{fig:acc_eff}(\subref{fig:eff_Peer})-(\subref{fig:eff_modRosPeer}),
the RosPeer schemes and its modified versions surpass the already existing LRSIF
versions of the implicit integration schemes for DREs. Further, it is noteworthy
that the fourth-order peer schemes do not reach better error levels that was
already visible from
Figures~\ref{fig:acc_eff}(\subref{fig:acc_Peer})-(\subref{fig:acc_modRosPeer}). A
more detailed investigation of all methods up to order \(2\) and in particular
the peer schemes can be found in~\cite{Lan17}.

\subsection{Convection-diffusion - Small-Scale LTV}
\label{sec:conv-diff-small}
\begin{figure}[t]
  \centering
  \begin{subfigure}{.325\linewidth}%
    \includegraphics{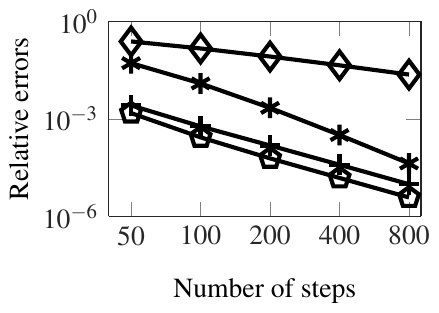}
    \caption{Peer(1-4)}\label{fig:FDM81_acc_Peer}
  \end{subfigure}\hskip-.25em%
  \begin{subfigure}{.325\linewidth}\hskip-.5em%
    \includegraphics{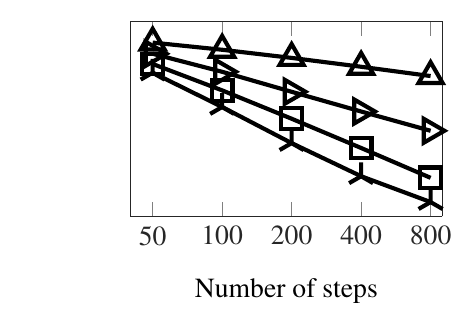}
    \caption{RosPeer(1-4)}\label{fig:FDM81_acc_RosPeer}
  \end{subfigure}%
  \begin{subfigure}{.325\linewidth}\hskip-.75em%
    \includegraphics{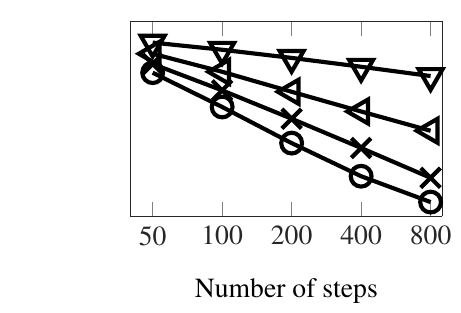}
    \caption{mRosPeer(1-4)}\label{fig:FDM81_acc_modRosPeer}
  \end{subfigure}\vskip.25em%

  \begin{subfigure}{.325\linewidth}
    \includegraphics{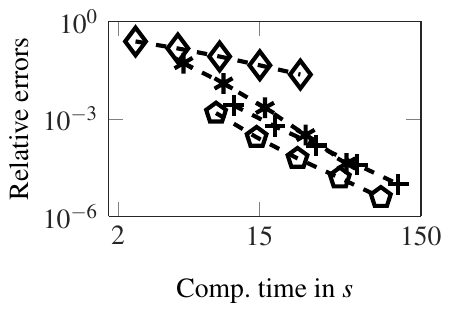}
    \caption{Peer(1-4)}\label{fig:FDM81_eff_Peer}
  \end{subfigure}\hskip-.25em%
  \begin{subfigure}{.325\linewidth}\hskip-.5em%
    \includegraphics{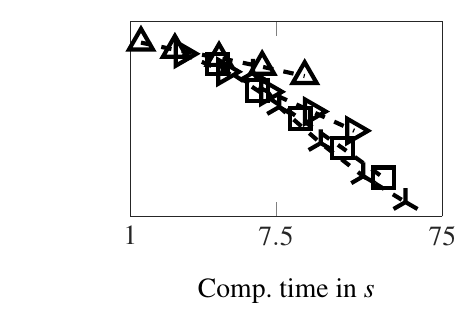}
    \caption{RosPeer(1-4)}\label{fig:FDM81_eff_RosPeer}
  \end{subfigure}%
  \begin{subfigure}{.325\linewidth}\hskip-.75em%
    \includegraphics{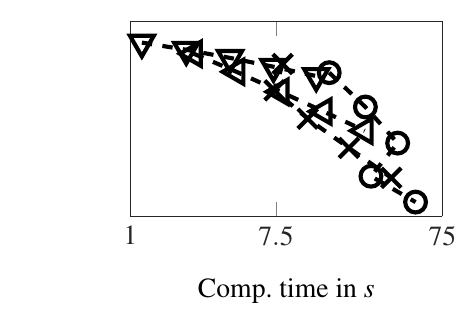}
    \caption{mRosPeer(1-4)}\label{fig:FDM81_eff_modRosPeer}
  \end{subfigure}%
  \medskip
  \begin{center}\vskip-1em
    \includegraphics[scale=1.25]{legend}
  \end{center}\vskip-1em
  \caption{Convection-diffusion LTV: Accuracy and efficiency plots}
  \label{fig:FDM81_acc_eff}
\end{figure}
The second example is a convection-diffusion model problem originating from a
centered finite differences discretization of the partial differential equation
\begin{align}
  \dot{v}=-\Delta v-f_{1}\frac{\partial v}{\partial\xi_{1}}-f_{2}
  \frac{\partial v}{\partial\xi_{2}}-f_{3}=0,\label{eq:conv-diff}
\end{align}
for \(v=v(\xi_{1},\xi_{2})\) defined on the unit square \(\Omega=(0, 1)^{2}\)
with homogeneous Dirichlet boundary conditions. Here, \(f_{i},~i=1,2,3\), are
functions depending on \(\xi_{1},\xi_{2}\) and are often referred to as
convection and reaction terms. The system matrices \(A\) and \(B,C\) are
generated by the \matlab ~routines \texttt{fdm\_2d\_matrix} and
\texttt{fdm\_2d\_vector}, respectively, from \lyapack~\cite{Pen00b} with
\(n_{0}=9\) equidistant grid points for each spatial dimension, resulting in
\(n=n_{0}^{2}=81\) unknowns, and the convection and reaction terms are chosen as
\(f_{1}=20,~f_{2}=5,~f_{3}=0\). Further, the model represents a single-input
single-output (SISO) system with \(m=1\) input and \(q=1\) output. The regions,
where \(B\) and \(C\) act are restricted to the lower left corner
\(\xi_{1}\in (0,0.35),\xi_{2}\in (0,0.35)\) for the input and the upper area
defined by \(\xi_{1}\in (0,1),\xi_{2}\in (0.95,1)\) for the output,
respectively. In order to obtain an LTV model, we introduce an artificial
time-variability $\mu(t)=\frac{3}{4}\sin(8\pi t)+1\in[0.25,1.75]$ to the system
matrix $A$. As a result, we obtain a time-varying system with constant matrices
$E,B,C$ and a time dependent matrix $A(t)=\mu(t)A$. The model is simulated for
the time interval \([0,0.5]\ s\) with the time step sizes
\(\tau\in\{\frac{1}{100}, \frac{1}{200}, \frac{1}{400}, \frac{1}{800},
\frac{1}{1600}\}\ s\), resulting in \(\{50, 100, 200, 400, 800\}\) steps,
respectively. As for the previous example,
Figures~\ref{fig:FDM81_acc_eff}(\subref{fig:FDM81_acc_Peer})-(\subref{fig:FDM81_acc_modRosPeer})
show the error behavior with respect to the several time step sizes used. Here,
the predicted convergence behavior is clearly visible except for the Peer(4)
scheme. The efficiency plots are presented in
Figures~\ref{fig:FDM81_acc_eff}(\subref{fig:FDM81_eff_Peer})-(\subref{fig:FDM81_eff_modRosPeer}). For
this example, again the peer schemes show best performance with respect to the
achieved accuracy. Additionally considering the computational effort of the
integration schemes, the BDF methods show best performance up to order 3. See
also Table~\ref{tab:fdmLTV_timings}. For large-scale model problems, the
computational effort for solving the ARE inside the implicit schemes will become
more expensive compared to the ALE solves within the linear implicit time
integrators such that the latter will become more effective.

{\rowcolors{3}{white}{gray!25}
  \begin{table}[t]
    \centering
    \caption{Convection-diffusion LTV: Computational timings and relative errors
      with respect to the reference solution for $\tau=6.25e$-4 \(s\), 800
      steps.}\label{tab:fdmLTV_timings}
    \begin{tabular}{lrc}
      \rowcolor{gray!50}
      Method&Time in $s$&Rel. Frobenius err.\\
      BDF(1)&25.07&2.32e-02\\
      BDF(2)&23.33&6.79e-04\\
      BDF(3)&23.05&7.34e-05\\
      BDF(4)&23.08&2.91e-05\\
      Ros1&12.57&2.09e-02\\
      Ros2&48.50&2.87e-03\\
      Ros4&62.18&4.36e-04\\
      Mid&29.46&1.91e-04\\
      Trap&29.11&2.13e-04\\
      Peer(1)&26.92&2.32e-02\\
      Peer(2)&51.93&4.26e-05\\
      Peer(3)&84.82&3.84e-06\\
      Peer(4)&108.81&9.81e-06\\
      RosPeer(1)&11.16&2.09e-02\\
      RosPeer(2)&22.15&4.32e-04\\
      RosPeer(3)&33.28&1.54e-05\\
      RosPeer(4)&45.03&2.77e-06\\
      mRosPeer(1)&13.09&2.09e-02\\
      mRosPeer(2)&25.60&4.32e-04\\
      mRosPeer(3)&37.00&1.54e-05\\
      mRosPeer(4)&51.74&2.77e-06
    \end{tabular}
  \end{table}
}
\subsection{Convection-Diffusion - Large-Scale LTI}
\label{sec:convection-diffusion}
{\rowcolors{3}{white}{gray!25}
  \begin{table}[t]
    \centering
    \caption{Convection-diffusion LTI: Computational timings for $\tau=6.25e$-4,
      480 steps.}\label{tab:fdmLTI_timings}
    \begin{tabular}{lr}
      \rowcolor{gray!50}
      Method&Time in $s$\\
      BDF(1)&1\,260.43\\
      BDF(2)&1\,038.41\\
      BDF(3)&870.76\\
      BDF(4)&813.07\\
      Ros1&1\,107.56\\
      Ros2&5\,779.17\\
      Ros4&10\,571.82\\
      Mid&793.03\\
      Trap&796.44\\
      Peer(1)&1\,239.61\\
      Peer(2)&1\,322.49\\
      Peer(3)&2\,068.50\\
      Peer(4)&2\,652.84\\
      RosPeer(1)&583.14\\
      RosPeer(2)&561.36\\
      RosPeer(3)&740.52\\
      RosPeer(4)&913.28\\
      mRosPeer(1)&584.07\\
      mRosPeer(2)&543.07\\
      mRosPeer(3)&647.03\\
      mRosPeer(4)&887.04
    \end{tabular}
  \end{table} }

The third example is again the convection-diffusion model~\eqref{eq:conv-diff}
from Example 2. Here, the convection and reaction terms
\(f_{1}=50,~f_{2}=10,~f_{3}=0\) and no additional artificial time-variability
are used. Further, \(n_{0}=45\) grid nodes in each direction, yielding a system
dimension of \(n=2\,025\), are considered. The model is simulated for the time
interval \([0,0.3]\ s\) with time step sizes
\(\tau\in\{\frac{1}{100}, \frac{1}{200}, \frac{1}{400}, \frac{1}{800},
\frac{1}{1600}\}\ s\), and \(\{30, 60, 120, 240, 480\}\) steps,
respectively. Due to the system size, no reference solution is computed. Similar
to the previous examples, Table~\ref{tab:fdmLTV_timings} shows the
computational timings for the several integration schemes. Again, the
Rosenbrock-type peer schemes up to order 3 come up with the lowest computational
times. It can also be seen that for this autonomous SISO system, the
reformulated Rosenbrock-type schemes (mRosPeer) outperform their counterparts
given in the original formulation.

\section{Conclusion}\label{sec:conc}
In this contribution, the classes of implicit and Rosenbrock-type peer methods
have been applied to matrix-valued ODEs. Further, a reformulation of the latter
has been proposed in order to avoid a number of Jacobian applications to the
currently computed stage variables. An efficient low-rank formulation in terms
of the low-rank symmetric indefinite factorization (LRSIF) has been
presented. The performance of the peer methods was presented for three different
examples. It has been shown that the Rosenbrock-type schemes and their
reformulated version outperform their classical implicit one- and multi-step
opponents with respect to the relation of accuracy and computational effort in
most cases. Thus, the peer methods and in particular Rosenbrock-type schemes
make an important contribution to the efficient low-rank based solution of
differential Riccati equations and most probably differential matrix equations,
in general.

\section*{Acknowledgements}
\paragraph{Financial Support} This research was funded by the Deutsche
Forsch\-ungsgemeinschaft DFG in subproject A06 “Model Order Reduction for
Thermo-Elastic Assembly Group Models” of the Collaborative Research Center/
Transregio 96 “Thermo-energetic design of machine tools – A systemic approach to
solve the conflict between power efficiency, accuracy and productivity
demonstrated at the example of machining production”.

\paragraph{Special Thanks} goes to
Prof. R. Weiner\footnote{\url{http://www.mathematik.uni-halle.de/wissenschaftliches_rechnen/ruediger_weiner/}}
and his group at the Martin-Luther-Universität Halle for helpfull discussions on
the peer methods and in particular for providing the coefficients for numerous
implicit peer schemes.
\bibliographystyle{siamplain}
\bibliography{refs}

\end{document}